\documentclass{article}

\usepackage{amssymb,dsfont,latexsym}

\newtheorem{thm}{Theorem}

\newtheorem{prop}[thm]{Proposition}
\newtheorem{lem}[thm]{Lemma}
\newtheorem{cor}[thm]{Corollary}

\newcommand\enu[1]{\smallskip\newline\makebox[5mm][l]{\rm(#1)}}

\newcommand\bp{\noindent{\it Proof.}\ }

\begin{document}

\author{Erling St{\o}rmer}

\date{11-10-2007 }

\title{Separable states and positive maps}

\maketitle
\begin{abstract}
Using the natural duality between linear functionals on tensor
products of C*-algebras with the trace class operators on a
Hilbert space $H$ and linear maps of the C*-algebra into $B(H)$,
we study the relationship between separability, entanglement and
the Peres condition of states and positivity properties of the
linear maps.

\end{abstract}

\section*{Introduction}

 In an earlier paper \cite{st3} we studied the duality between
 linear functionals $\tilde\phi $ on a tensor product
 $A\widehat\otimes{\1T}$ of an operator system $A$ and the trace class
 operators $\1T$ on a Hilbert space $H$, and bounded linear maps
 $\phi\colon A \to B(H) $ given by the formula $\tilde\phi(a\otimes b)=
 Tr(\phi(a)b^t)$. The main emphasis was on positivity properties
 of $\tilde\phi$ on cones in $A\widehat\otimes{\1T}$ obtained by classes of
 positive maps. In the present paper we shall see how this study
 yields a natural framework for the study of separable states of
 $A\widehat\otimes{\1T}$, for example we recover results of Horodecki
 et.al \cite {3 Hor} and Horodecki, Shor and Ruskai \cite{HSR} on
 characterizations of separable states. In addition we shall obtain
 characterizations of states on $A\widehat\otimes{\1T}$ satisfying the
 Peres condition, viz $\rho\circ(\iota\otimes t)$ is positive,
 where t is the transpose map and $\iota$ the identity map. In
 particular we see that nondecomposable maps yield natural
 examples of entangled states which satisfy the Peres condition; for this see also \cite{Kye1},
 \cite{Kye2}.In the last section we study the definite set of a positive map
$\phi$ on a $C^*$-algebra $A$, i.e. the set of self-adjoint
operators in $A$ such that $\phi(a^2)=\phi(a)^2,$ and show that if
$\tilde\phi$ is a separable state, then the image of the definite
set is a $C^*$-subalgebra of the center of the $C^*$-algebra
generated by $\phi(A)$. As a corollary we obtain a decomposition
result for separable states in the finite dimensional case.

The author is indebted to E.Alfsen for his careful reading of the manuscript 
and several useful comments.

\section{Cones and states}

In this section we recall notation and concepts from \cite{st3}
and show a general characterization of separable states close to
that in \cite{HSR}. For more details on the following see
\cite{st3}.

By an \textit{operator system} we shall mean a norm closed
self-adjoint set $A$ of bounded operators on a Hilbert space
containing the identity.  We denote by $A \odot B(H)$ the
algebraic tensor product of $A$ and $B(H)$  and by $A\otimes B(H)$
the closure of $A \odot B(H)$ in the operator norm. If $\1T$
denotes the trace class operators on $H$, then
$A\widehat\otimes{\1T}$ denotes the projective tensor product of
$A$ and $\1T$.  We denote by $B(A,H),$ (resp. $B(A,H)^+$ ) the
bounded (resp. positive) linear maps of $A$ into $B(H)$. The
\textit{BW-topology} on $B(A,H)$ is the topology of bounded
pointwise weak convergence, i.e. a net $(\phi_\nu)$ converges to
$\phi$ if it is uniformly bounded, and $\phi_\nu(a) \rightarrow
\phi(a)$ weakly for all $a\in A$.  We denote by $t$  the transpose
map on $B(H)$ with respect to some orthonormal basis for $H$. Then
by abuse of notation we get that the transpose map on $B(K)\otimes
B(H)$ is $t\otimes t$. We shall also denote by $Tr$ the usual
trace on $B(H)$ which takes the value 1 on minimal projections. We
recall Lemma 2.1 in \cite{st3}.

\begin{lem}\label{lem 2.1}
With the above notation there is an isometric isomorphism $\phi
\rightarrow\tilde\phi $ between $B(A,H)$ and
$(A\widehat\otimes{\1T})^*$ given by
$$
\tilde\phi(a\otimes b)= Tr(\phi(a)b^t), \ a\in A, b\in \1T.
$$
Furthermore, $\phi \in B(A,H)^+$ if and only if $\tilde\phi$ is
positive on the cone $A^+\widehat\otimes{\1T}^+$ generated by
operators of the form $a\otimes b$ with $a$ and $b$ positive.
\end{lem}

We recall Definition 2.3 in \cite{st3}.  It says that a BW-closed
subcone $K\neq 0$ of $B(B(H),H)^+$ is a \textit{mapping cone} if
it has a BW-dense subset of ultra weakly continuous maps and is
invariant in the sense that if $\alpha\in K$, and $a,b \in B(H)$
then the map $x\rightarrow a\alpha(bxb^*)a^*$ belongs to $K.$
Three mapping cones will be of special interest in the
following, namely $B(B(H),H)^+, CP(H)$ - the set of completely
positive maps in $B(B(H),H)$, and $S(H)$ - the BW-closed cone
generated by maps of the form
$$
x\rightarrow \sum_{i=1}^{n}\omega_i(x)a_i,
$$
where $\omega_i$ is a normal state on $B(H)$ and $a_i\in B(H)^+$.
The latter maps are said to be of "Holevo form" in \cite{HSR}. By
Lemma 2.4 in \cite{st3} $S(H)$ is the minimal mapping cone and
$B(B(H),H)^+$ the maximal one.

If $K$ is a mapping cone and $A$ an operator system as before, we
denote by $P(A,K)$ the cone
$$
P(A,K) = \{x\in A\widehat\otimes{\1T}: \iota\otimes\alpha(x)\geq
0\  \forall \alpha\in K\}.
$$
By Lemma 2.8 in \cite{st3} $P(A,K)$ is a proper norm closed convex
cone in $A\widehat\otimes{\1T}$ containing the cone
$A^+\widehat\otimes{\1T}^+$. A map $\phi \in B(A,H)$ is said to be
\textit{K-positive} if
$$
\tilde\phi(\sum a_i\otimes b_i)=\sum Tr(\phi(a_i)b_i^t) \geq 0 \
whenever \ \sum a_i\otimes b_i \in P(A,K).
$$
By Theorem 3.2 in \cite{st3} $\phi$ is completely positive if and
only if $\tilde\phi$ is positive on the cone
$(A\widehat\otimes{\1T})^+$, the closure of the positive operators
in $A\odot\1T$,  if and only if $\phi$ is $CP(H)-$
positive.

If $C\subseteq V$ and $D\subseteq W$ are closed convex cones in
two real locally convex vector spaces in duality, we denote by
$C^*$ (resp. $D^*$) the set of $w\in W$ such that $<v,w>\ \geq 0
 \; \forall v\in C$, (resp. $v\in V$ such that $<v,w>\ \geq 0 \ \forall
w\in D)$. Thus $\phi$ is K-positive if and only if $\tilde\phi \in
P(A,K)^*$. By a straightforward application of the Hahn-Banach
Theorem for closed convex cones, see e.g. \cite{AS}, Prop. 1.32,
we have
$$
P(A,K) = P(A,K)^{**}.
$$
We say a positive linear functional $\rho$ on $A\otimes B(H)$ is \textit{separable} if
it belongs to the norm closure of positive sums of states of the
form $ \sigma\otimes\omega$, where $\sigma$ is a state of $A$ and
$\omega$ a normal state of $B(H).$ Otherwise $\rho$ is called
\textit{entangled}. We denote the set of separable states by
$S(A,H).$ It is a norm closed cone in $
(A\widehat\otimes{\1T})^*$. As for $P$ above $S(A,H) = S(A,H)^{**}.$Our next result is closely related to
Theorem 2 in \cite{HSR}.

\begin{thm}\label{TH2}
Let $A$ be an operator system and $\phi\in B(A,H).$ Then the
following conditions are equivalent:
 \enu{i} $\tilde\phi$ is a
separable positive linear functional.
 \enu{ii} $\phi$ is
$S(H)-$positive.
 \enu{iii} $\phi$ is a BW-limit of maps of the
form $x\rightarrow \sum_{i=1}^n \omega_i(x)b_i$ with $\omega_i$ a
state of $A$, and $b_i\in B(H)^+$.
\end{thm}
\bp $(i)\Leftrightarrow(ii)$. Let $Sn$ denote the positive normal
linear functionals on $B(H)$, and let $x=\sum x_i\otimes y_i \in
A\odot B(H).$ Then
$$
x\in P(A,S(H))
$$
$$ \Leftrightarrow (\iota\otimes b\omega)(x) \geq 0 \
\forall \omega\in Sn, b\geq 0
$$
$$
\Leftrightarrow \sum x_i\omega(y_i)\otimes b = \sum
x_i\otimes\omega(y_i)b\geq 0 \; \forall \omega\in Sn, b\geq 0
$$
$$
\Leftrightarrow\sum x_i\omega(y_i)\geq 0  \; \forall \omega\in Sn
$$
$$
\Leftrightarrow\rho\otimes\omega(x)=\sum
\rho(x_i)\omega(y_i)=\rho(\sum x_i\omega(y_i))\geq 0 \; \forall
\omega\in Sn, \rho\in A^{*+}
$$
$$
\Leftrightarrow \eta(x)\geq 0 \; \forall \eta\in S(A,H)
$$
$$
\Leftrightarrow x\in S(A,H)^*.
$$
Thus $\phi$ is $S(H)$-positive if and only if $\tilde\phi \in
P(A,S(H))^* = S(A,H)^{**} = S(A,H))$, proving that $(i)\Leftrightarrow (ii)$. The
equivalence $(ii) \Leftrightarrow (iii)$ follows from Theorem 3.6
in \cite{st3}, since a map $\alpha\in S(H)$ if and only if $t
\circ \alpha \circ t \in S(H)$. The proof is complete.
\vskip0.5cm
In \cite{HSR} maps like $x\rightarrow\sum \omega_i(x)b_i$ are
called "entanglement breaking".

 It is possible to give a direct
proof of a less general form of the equivalence
$(i)\Leftrightarrow (iii)$ above. Suppose $\phi(a)=\sum
\omega_i(a)b_i$ for $a\in A, b_i\in B(H)^+, \omega_i$ state of
$A$. Then
$$
\tilde\phi(a\otimes b)=Tr(\phi(a)b^t)=\sum Tr(\omega_i(a)b_i
b^t)=\sum \omega_i(a) Tr(b_i b^t)=\sum \omega_i(a)\rho_i(b),
$$
where $\rho_i(b)=Tr(b_i b^t)$ is a positive linear functional.
Thus $\tilde\phi$ is separable.

Conversely, if $\tilde\phi = \sum \omega_i\otimes\rho_i$, let
$\tilde\rho_i (b)=\rho_i(b^t)= Tr(b_i b)$. Then we have
$$
Tr(\phi(a)b^t)=\tilde\phi(a\otimes b)=\sum
\omega_i(a)\rho_i(b)=\sum \omega_i(a)\tilde\rho_i(b^t)=\sum
Tr(\omega_i(a)b_i b^t).
$$
This holds for all $b\in \1T$, hence $\phi(a)=\sum
\omega_i(a)b_i.$

\begin{cor}\label{cor2.4}
Let $H$ be separable and $\phi\in B(A,H)^+.$ Suppose $\phi(A)$ is
contained in an abelian $C^*-$algebra. Then $\tilde\phi$ is
separable.
\end{cor}
\bp By hypothesis there is an abelian von Neumann algebra
$B\subseteq B(H)$ such that $\phi\colon A\to B.$ Let $(B_n)$ be an
increasing sequence of finite dimensional von Neumann subalgebras
of $B$ such that $\bigcup_n B_n$ is weakly dense in $B.$ Let
$E_n\colon B\to B_n$ be normal conditional expectations such that
$E_{n-1}|_{B_n}\circ E_n = E_{n-1}.$ Then $\phi(x)=wlim_n E_n\circ
\phi(x)$ for all $x\in A.$ Since $B_n$ is finite dimensional,
$E_n\circ\phi(x)=\sum \omega_i^n(x)e_i^n,$ where $\omega_i^n$ are
positive linear functionals  on $A$ and $e_i^n$ are minimal
projections in $B_n.$ Since $\phi$ is a BW-limit of the sequence
$E_n\circ\phi$, $\tilde\phi$ is separable by Theorem \ref{TH2}.
The proof is complete.
\vskip0.5cm
A celebrated necessary condition for a state $\rho$ on
$A\widehat{\otimes} \1T$ to be separable is the \textit{Peres
condition}, i.e. $\rho\circ(\iota\otimes t)\geq 0.$ A map $\phi\in
B(A,H)$ is said to be \textit{copositive} if $t\circ\phi$ is
completely positive.

\begin{prop}\label{pro4}
Let $\phi\in B(A,H).$ Then $\tilde\phi$ satisfies the Peres condition
if and only if $\phi$ is both completely positive and copositive.
\end{prop}
\bp If $a\in A$ and $b\in\1T$ we have, since the trace is
invariant under transposition,
$$
\tilde\phi(a\otimes b^t)
=Tr(\phi(a)b)=Tr(t\circ\phi(a)b^t)=(t\tilde\circ\phi)(a\otimes b).
$$
Thus $\tilde\phi$ satisfies the Peres condition if and only if both
$t\tilde\circ\phi$ and $\tilde\phi$ are positive.  Using Theorem 3.2 in \cite{st3} 
this holds  if and only if
$t\circ\phi$ and $\phi$ are completely positive, hence
if and only if $\phi$ is both completely positive
and copositive.

\section {States on $B(K)\otimes B(H)$}.

In this section we study the case when the operator system $A$ equals $B(K)$ for
a Hilbert space $K$. But first we consider the finite dimensional case.
Let $M_n=M_n(\0C)$ denote the complex $n\times n$ matrices, and let
$\phi\colon M_n\to M_m$, so $\phi\in B(A,\0C^m),$ where $A=M_n$
and $H=\0C^m$. Let $(e_{ij})$ be a complete set of matrix units in
$M_n.$ Then the \textit{Choi matrix} for $\phi$ is
$$
C_{\phi} =\sum e_{ij}\otimes \phi(e_{ij})=\iota\otimes\phi(P)\in
M_n\otimes M_m,
$$
where $\frac{1}{n}P$ is the 1-dimensional projection
$\frac{1}{n}\sum e_{ij}\otimes e_{ij},$ - the so-called maximally
entangled state, see \cite{Cho1}. Denote by $\phi^t$ the map
$t\circ\phi\circ t$, where $t$ denotes the transpose map in either
$M_n$ or in $M_m.$ Then $\phi$ is completely positive if and only
if $\phi^t$ is completely positive.  It was shown by Choi
\cite{Cho1} that $\phi$ is completely positive if and only if
$C_{\phi}$ is positive.  We use the convention that the density matrix
for a state $\rho$  is the positive matrix $h$ such that $\rho(x)=Tr(hx)$.

\begin{lem}\label{lem3.1}
$C_{\phi^t}$ is the density matrix for $\tilde\phi.$
\end{lem}

\bp Let $a\in M_n, b\in M_m.$ Since the transpose $t$ on $M_n
\otimes M_m$ is the tensor product of the transpose operators on
$M_n$ and $M_m$, we have
\begin {eqnarray*}
Tr(C_{\phi^t}a\otimes b)
&=& \sum Tr(e_{ij}\otimes\phi^t(e_{ij})(a\otimes b))\\
&=& \sum Tr(e_{ji}\otimes \phi(e_{ij}^t)(a^t\otimes b^t))\\
&=& \sum Tr(e_{ji}a^t) Tr(\phi(e_{ji})b^t)\\
&=& \sum a_{ji}Tr(e_{ji}\phi^* (b^t))\\
&=& Tr(a\phi^* (b^t))\\
&=& \tilde\phi(a\otimes b).
\end{eqnarray*}
In the above computation $\phi^*$ is the adjoint of $\phi$ as an
operator between $M_n$ and {$M_m$ considered as the
trace class operators on $\0C^n$ and $\0C^m$ respectively. The
proof is complete.

\begin{lem}\label{lem3.2}
Let $H=\0C^m$ and $\phi\in B(M_n,H).$ Then $\phi$ is positive if
and only if $C_{\phi^t}\in P(M_n,S(H))$, if and only if
$C_{\phi}\in P(M_n,S(H))$. Hence $P(M_n,S(H))= \{C_\phi : \phi
\geq 0\}.$
\end{lem}
\bp By Theorem \ref{TH2}, or rather the proof of the equivalence
$(i)\Leftrightarrow (ii),$
$$
C_{\phi^t}\in P(M_n,S(H))=S(M_{n},H)^*
$$
$$
\Leftrightarrow Tr(C_{\phi^t} a\otimes b)\geq 0 \ \forall a\in
M_n^+, b\in M_m^+
$$
$$
\Leftrightarrow \phi \geq 0
$$
by Lemma \ref{lem 2.1}, proving the first statement. Since
$\phi\geq 0$ if and only if $\phi^t \geq 0,$ the above is
equivalent to $C_\phi$ being in $P(M_n,S(H))$.

Each element $x\in P(M_n,S(H))$ defines a linear functional $\rho$
on $M_n\otimes M_m$ by $\rho(y)=Tr(xy).$ By Lemma \ref{lem 2.1}
there is $\phi \in B(M_n,\0C^m)$ such that $\rho(a\otimes
b)=Tr(\phi(a)b^t),$ hence by Lemma \ref{lem3.1} and the first part of 
the proof, $x=C_{\phi^t}$with $\phi \geq 0.$ Thus the last statement follows, 
completing the proof.
\vskip0.5cm
We shall now apply the finite dimensional results to study states on $B(K)\otimes B(H)$
and to prove an infinite dimensional extension of the Horodecki Theorem \cite{3 Hor}.
Recall that a state and a positive linear map on a Von Neumann algebra are said to be
normal if they are weakly continuous on bounded sets.

\begin{thm}\label{thm3.3}
Let $\rho$ be a normal state on $B(K)\otimes B(H)$ with $K$ and $H$ Hilbert spaces
and with density operator $h$.Then $\rho$ is separable if and only if
 $\iota\otimes \psi(h)\geq 0$ for all normal positive maps
 $\psi\colon B(H)\to B(K).$
\end{thm}
\bp Suppose $\rho$ is separable and normal. Then $\rho\circ ( \iota\otimes \phi)$ 
is a  normal state for all unital normal positive maps $\phi\colon B(K) \to B(H).$
Let $\psi$ be as in the statement of the theorem. Then the adjoint map $\psi^{*}$ is
 a positive map of the trace class operators on $K$ into those on $H.$
Thus if $x\geq 0$ is of finite rank in $B(K\otimes K) = B(K)\otimes B(K)$, then
$$
Tr((\iota\otimes\psi)(h)x) =   Tr(h(\iota\otimes\psi^{*})(x)) = \rho(\iota\otimes\psi^{*}(x))
\geq 0,
$$
hence $\iota\otimes\psi(h)\geq 0.$ 

To show the converse we first assume $K$ and $H$ are finite dimensional.  Then
by Lemma \ref{lem3.2} $P(M_n,S(H))= \{C_\phi : \phi \geq 0\}.$
Thus by Theorem \ref{TH2} and Lemma \ref{lem3.1} $\rho$ is
separable if and only if for all positive $\phi\colon B(K)\to B(H)$
$$
Tr((\iota\otimes\phi^*)(h)P) = Tr(h(\iota\otimes\phi)(P))=Tr(hC_\phi)\geq 0 ,
$$
where $P$ is the rank one matrix such that $C_\phi=\iota\otimes\phi(P).$ 
Since $P\geq 0,$ and by assumption $\iota\otimes\phi^{*}(h) \geq 0$,
it follows that $\rho$ is separable.  

We next consider the general case when $K$ and $H$ may be infinite dimensional. 
Assume $\iota\otimes\psi(h)\geq 0$ for all normal  $\psi\colon B(H)\to B(K).$
Since the maps $\psi_{f}(x)=\psi(fxf)$ are positive for all finite dimensional projections $f$, it i clear that
 $\iota\otimes\psi((e\otimes f) h (e\otimes f))\geq 0$ for all normal positive maps
$\psi\colon B(H) \to B(K)$ with $e$ a finite dimensional projection in $B(K))$.    Let
$$
\psi_{e\otimes f}(y) = e\psi(fyf)e,  \ y\in B(H).
$$
Then $\iota \otimes \psi_{e\otimes f}(h)\geq 0.$ Now every normal positive map
$\phi\colon B(fH) \to B(eK)$ is of the form $ \psi_{e\otimes f}$ with $\psi$ as
above, because we can define $\phi\colon B(H) \to B(K)$ by $\psi(x) = \phi(fxf).$
Thus by the part of the proof on the finite dimensional case, the positive linear functional
$\omega(x) = \rho((e\otimes f) x( e\otimes f))$ is separable on $B(eK)\otimes B(fH)$.
Since this holds for all finite dimensional projections $e$ and $f$ and $\rho$ is normal, it
follows that $\rho$ is separable.  The proof is complete.
\vskip 0.5cm

We expect that the above theorem can be generalized to Von Neumann algebras other than
$B(K).$  If $A$ is a semi-finite Von Neumann algebra then so is $A\otimes B(H),$ hence 
each normal state on  $A\otimes B(H)$ has a density operator with respect to a trace, and 
the formulation of the theorem has a natural generalization.  In the type III case a formulation 
in terms of modular theory ought to be possible. 

\vskip 0.5 cm

  We next restate the Peres condition in terms of the density matrix of
the normal state $\rho.$

\begin{thm}\label{thm3.4}
Let $\rho$ be a normal state on $ B(K)\otimes B(H)$  with density operator $h,$ and let $t$ 
denote the transpose map of either $B(K)$ or $B(H)$.  Then the following conditions are
 equivalent: 
 \enu{i} $\rho$ satisfies the Peres condition.
 \enu{ii}$ \iota\otimes t(h)\geq 0$.
 \enu{iii} $t\otimes\iota (h)\geq 0$.
 \enu{iv} $h\in P(B(K),CP(H)) \bigcap P(B(K),copos(H))$,
 where $copos(H)$ denotes the copositive maps of $B(H)$ into
 itself.
 \end{thm}
\bp Assume (i).  Since the trace on $ B(K)\otimes B(H)$ 
is invariant under $\iota\otimes t$ , we have
$$
\rho\circ (\iota\otimes t)(a\otimes b) = Tr(h(\iota\otimes t)(a\otimes b))
= Tr(\iota\otimes t(h)(a\otimes b)).
$$
Since $\rho\circ \iota\otimes t\geq 0$ it follows that $\iota\otimes t(h)\geq 0$.

Conversely, if (ii) holds then by the above computation $\rho\circ (\iota\otimes t)\geq 0,$
hence (i) holds.The equivalence of (ii) and (iii) follows since
 $t\otimes\iota(h) = t\otimes t (\iota\otimes t (h)),$
and the fact that $t\otimes t$ is an order-automorphism.

We have
$$
P(B(K),copos(H)) = \{x\in B(K)\otimes B(H): \iota \otimes \phi(x)\geq 0 \ \forall \
copositive \ \phi\}
$$ 
$$
=\{x\in B(K)\otimes B(H): \iota \otimes t(x)\geq 0\},
 $$
 because a copositive map is the composition of a completely positive map and the 
 transpose map.  Thus (ii) is equivalent to (iv), completing the proof.
\vskip0.5cm

Let $A$ be a $C^*$ - algebra. Then  a  map $\phi\in B(A,H)$ is called \textit{decomposable} if it is
the sum of a completely positive map and a copositive map.
Otherwise $\phi$ is called \textit{nondecomposable}.  Since a map
$\phi\in B(A,\0C^n)$ is completely positive if and only if
$\iota\otimes\phi\colon M_{n}\otimes A \to M_n\otimes M_n$ is
positive \cite{ER},Lemma 5.1.3, it follows from \cite{st2} that
$\phi\in B(A,\0C^n)$ is decomposable if and only if whenever $h$
and $t\otimes \iota (h)$ belong to $(M_{n}\otimes A)^+$ then
$\iota\otimes\phi(h)\geq 0.$  Thus $\phi$ is nondecomposable if and
only if there exists $h\in (M_n\otimes A)^+$ such that
$t\otimes\iota(h)\geq 0$ while $\iota\otimes\phi(h)$ is not
positive. Suppose that $A=B(H), \phi$ normal, and $h$ as above.
Then there exists by normality of $\phi$ a finite dimensional projection
$f\in B(H)$ such that $\iota\otimes \phi((1\otimes f)h(1\otimes f)$ is
not positive.  We can thus assume $h$ is of finite rank.  Normalizing $h$ we thus have by 
Theorem \ref{thm3.4} that the state $\rho(x)=Tr(hx)$ satisfies the Peres condition, while by
Theorem \ref{thm3.3} $\rho$ is entangled.  We have thus proved

\begin{thm}\label{thm3.5}
Let $\phi\colon B(H)\to M_n$ be normal positive and nondecomposable. 
Then there exists a normal state $\rho$ on $B(H)\otimes M_n $ with density operator $h$ such
that $t\otimes\iota (h)\geq 0$, while $\iota\otimes\phi(h)$ is not
positive. Hence $\rho$ is entangled but satisfies the Peres condition.
\end{thm}

An explicit example of the situation in the above theorem is given
in \cite{st2} and \cite{Cho3}. Then $ dim H =n=3$, and $\phi\colon M_3\to M_3$ is the
nondecomposable Choi map \cite{Cho2}. Other examples can be found in 
\cite{Kye1} and \cite{Kye2}.  A large class of
nondecomposable maps are the projections onto spin factors of
dimension greater than 6, or more generally, positive projections
onto nonreversible Jordan algebras, see \cite{st1}.  See \cite{T} for
another class of nondecomposable maps. Another result close to the above 
theorem can be found in \cite{BFP}. Previous
examples of entangled states which satisfy the Peres condition
have been exhibited by P.Horodecki \cite{PHor}.

\section{Definite sets}

If $A$ and $B$ are C*-algebras, and $\phi\colon A\to B$ is a
positive map of norm $\leq 1$ then the (self-adjoint)
\textit{definite set} $D_\phi$ of $\phi$ is the set of
self-adjoint operators in $A$ such that $\phi(a^2)=\phi(a)^2.$ If
$a\in D_\phi$ and $b\in A$ then
$\phi(ab+ba)=\phi(a)\phi(b)+\phi(b)\phi(a)$ and
$\phi(aba)=\phi(a)\phi(b)\phi(a),$ see \cite{st1}. We show in the
present section that if $\phi$ is of the form $\phi(x)=\sum
\omega_i(x)b_i$ as in Theorem \ref{TH2}, then $\phi(D_\phi)$ is
contained in the center of the C*-algebra generated by $\phi(A).$
In particular, if $\phi$ is faithful, then $D_\phi$ is abelian. As
a consequence we get a decomposition result for separable states.

\begin{thm}\label{thm4.1}
Let $A$ be a unital C*-algebra and $\phi\in B(A,H)^+$ with $\phi(1)=
1.$ Suppose $\phi$ is of the form $\phi(x)=\sum_{i=1}^n
\omega_i(x)b_i$ with $\omega_i$ states of $A$ and $b_i\in B(H)^+.$
Let $e$ be a projection in the definite set $D_\phi$ of $\phi,$
and put $f=1-e.$ Then $\phi(e)$ and $\phi(f)$ are projections in
$B(H)$ and satisfy
$$
\phi(x)=\phi(exe)+\phi(fxf)=\phi(e)\phi(x)\phi(e)+\phi(f)\phi(x)\phi(f)
$$
for all $x\in A.$ Hence $\phi(D_\phi)$ is an abelian C*-algebra
contained in the center of the von Neumann algebra generated by
$\phi(A).$ In particular, if $\phi$ is faithful then $D_\phi$ is
an abelian C*-algebra.
\end {thm}
\bp Since $e\in D_\phi, \phi(e)$ and $\phi(f)$ are mutually orthogonal  projections. Thus
$$
0=Tr(\phi(e)\phi(f))=Tr(\sum \omega_i(e)b_i \omega_j(f)b_j) = \sum
\omega_i(e)\omega_j(f)Tr(b_i b_j).
$$
Since each summand is positive we have
$$
\omega_i(e)\omega_j(f)Tr(b_i b_j)=0 \ \forall i,j.
$$
In particular
$$
\omega_i(e)\omega_i(f)Tr(b_{i}^2)=0 \ \forall i.
$$
Since $b_i\neq 0$ either $\omega_i(e)=0$ or $\omega_i(f)=0$ for
all i. In particular, $e$ or $f$ belongs to the left and right kernel
of $\omega_i,$ hence $\omega_i (exf)=\omega_i (fxe)=0$ for all $x$.
Thus $\omega_i(x)=\omega_i(exe)+\omega_i(fxf)$ for all $x$, so that
$$
\phi(x)=\phi(exe)+\phi(fxf)=\phi(e)\phi(x)\phi(e)+\phi(f)\phi(x)\phi(f),
$$
where the last equality follows since $e,f\in D_\phi.$

To show the last statement in the theorem we consider the
ultra-weakly continuous extension $\phi^{**}$ of $\phi$ to the
second dual $A^{**}$ of $A$. If $a\in D_\phi$ the abelian von
Neumann algebra generated by $a$ in $A^{**}$ is contained in
$D_{\phi^{**}}$ and is generated by its projections.  It thus
suffices to show that for each projection $e\in D_\phi, \phi(e)$
belongs to the commutant of $\phi(A)$. But this is immediate from
the above equation.

If $\phi$ is faithful then the restriction of $\phi$ to $D_\phi$
is an isomorphism, hence is abelian, since $\phi(D_\phi)$ is
abelian. The proof is complete.
\begin{cor}\label{cor4.2}
Let $A\subseteq B \subseteq B(H)$ be unital C*-algebras with $H$
separable. Suppose $\phi\colon B\to A$ is a conditional
expectation. Then $\tilde\phi$ is separable if and only if $A$ is
abelian.
\end{cor}
 \bp By Corollary \ref{cor2.4} if $A$ is abelian then $\tilde\phi$ is
 separable. Since $\phi$ is a conditional expectation, the self-adjoint part of $A$ equals
 the definite set $D_\phi$, hence by Theorem \ref{thm4.1} $A$ is
 abelian if $\tilde\phi$ is separable, completing the proof.
\vskip0.5cm
 Let $\tilde\phi=\sum \lambda_i\omega_i\otimes\rho_i$ be a
 faithful separable state on $M_n\otimes M_m$, which is a convex sum of
 states $\omega_i$ on $M_n$ and $\rho_i$ on $M_m$. By symmetry in $M_{n}$ and $M_{m}$
 in Lemma \ref{lem 2.1}, there exists a completely positive map $\psi\colon M_m\to
 M_n$ such that $\tilde\phi(a\otimes b)=Tr(a^t\psi(b)).$ Then by
 Theorem \ref{thm4.1} and the faithfulness of $\tilde\phi,D_\phi$ and $D_\psi$ are abelian
 C*-algebras. Let $(e_j)_{j=1,...,p}$ be minimal projections in
 $D_\phi$ and $(f_k)_{k=1,...,q}$ be minimal projections in
 $D_\psi.$ From the proof of Theorem \ref{thm4.1} the values of
 $\omega_i(e_j)$ and $\rho_i(f_k)$ are 0 or 1. In particular, the
 supports of $\omega_i$ and $\rho_i$ are contained in some $e_j$
 and $f_k$ respectively. Hence $e_j\otimes f_k$ are mutually
 orthogonal projections with sum 1 such that
 $$
 \tilde\phi(x)=\sum_{i,j} \tilde\phi(e_j\otimes f_k x e_j\otimes
 f_k),
 $$
 for all $x\in M_{n}\otimes M_{m}$

 We say $\tilde\phi$ is \textit{irreducible} if $D_\phi=D_\psi
 =\0R$ when we have cut down by the support of $\tilde\phi$,
 and we say a family $(\eta_i)$ of states are \textit{orthogonal} if their
 supports are mutually orthogonal. Summing up we have shown

 \begin{cor}\label{cor4.3}
 Every separable state on $M_{n}\otimes M_{m}$ is a convex sum  of orthogonal irreducible separable
 states.
 \end{cor}

Department of Mathematics, University of Oslo, 0316 Oslo, Norway.

e-mail: erlings@math.uio.no

\end{document}